\documentclass{article}

\usepackage[utf8]{inputenc}
\usepackage[T1]{fontenc}
\usepackage{lmodern}

\usepackage[a4paper, margin=1in]{geometry}

\usepackage{authblk}

\usepackage{mathtools}
\usepackage{amsmath, amsfonts, amssymb}

\usepackage{hyperref}
\usepackage{graphicx}

\usepackage{algorithm}
\usepackage{algpseudocode}

\usepackage{color}

\newcommand{\mut}{ \frac{\rho k}{\omega}}
\renewcommand{\div}{\mathrm{div}}
\newcommand{\T}{\intercal}
\newcommand{\R}{\mathbb{R}}

\title{Application of \texorpdfstring{$p$}{\textit{p}}-Laplacian relaxed steepest descent to shape optimization in two-phase flows}

\author[1]{Peter Marvin M{\"u}ller\footnote{\href{mailto:peter.marvin.mueller@tuhh.de}{peter.marvin.mueller@tuhh.de}}}
\author[2]{Martin Siebenborn}
\author[1]{Thomas Rung}

\affil[1]{Hamburg University of Technology}
\affil[2]{Hamburg University}

\begin{document}
	
\maketitle

\begin{abstract}
The paper is concerned with the minimal drag problem in shape optimization of merchant ships exposed to turbulent two-phase flows. Attention is directed to the solution of Reynolds Averaged Navier-Stokes equations using a Finite Volume method. Central aspects are the use of a  $p$-Laplacian relaxed steepest descent direction and the introduction of crucial technical constraints to the optimization procedure, i.e. the center of buoyancy and the displacement of the underwater hull. The example included refers to the frequently investigated Kriso container ship (KCS).    
\end{abstract}

\section{Introduction}
\label{sc:Introduction}
In this paper we formulate the minimal drag problem applied to the fluid dynamic shape optimization of merchant ships exposed to turbulent two-phase flows. Such optimization problems are additionally restricted by geometric constraints on the displacement and the center of gravity location.  
To this end, we propose an algorithm for handling these constraints based on first-order descent methods for this optimization problem as well as for the related sub-problems.

The shape optimization problem can be seen as an optimal control problem where the state is described by a set of partial differential equations (PDEs) which depend on a control.
The control, however, is given by the domain where again the state is defined on.
In addition the domain often has to fulfill additional geometric properties which yield a finite number of additional constraints.
For the constraint problem we apply the method of Lagrange multipliers where the derivative of the objective with respect to the control is expressed by primal (physical) and adjoint state variables.
It is well known that the adjoint approach is efficient for handling PDE constraint optimization problems featuring a large number of degrees of freedom (control variables) and can be realized with first order descent methods \cite{HPUU2009}.
Utilizing second-order methods, one could obtain an expression for updating the Lagrange multipliers in compliance with additional geometric constraints.
For example the method of mappings investigated in \cite{OS2020} allows to consider geometric constraints that depend on the state as the variables (primal and adjoint) are determined on the transformed domain.
A drawback of this approach, is that the whole optimality system has to be solved simultaneously, in line with -for example- one-shot methods \cite{OG2009}, which violates the sequential framework of many engineering CFD procedures \cite{FP2020}.
When solving the optimality system sequentially, the general approach is (1) compute the primal (physical) state, (2) compute the adjoint state variables or Lagrange multipliers and finally (3) determine a descent direction and shape deformation field, respectively.
The deformation field is obtained from the shape derivative of the objective function \cite{SZ1992, DZ2011}.

In \cite{MPRS2022} a first-order approach for the shape optimization problem with geometric constraints was investigated for geometric constraints that do not depend on the state and thus decouple from  the shape optimization problem. 
This allows exclusively consider the constraints when computing the shape deformation field by applying Newton's method with the Schur complement method for solving the related saddle point problem.
The fluid dynamic problems considered in the present study, however,  consider geometric constraints that also depend on the state. Strictly speaking, this does not allow to decouple the geometric constraints from the shape optimization problem.
Nevertheless, we will state that a decoupled strategy can be pursued for small step sizes  if the shape derivative is computed in the reference domain.

When considering the minimal drag problem of a free floating vessel it is necessary to conserve the displacement and center of buoyancy of the hull.
This is different to the typical geometric constraints that are given in an aerodynamic shape optimization problem, e.g. volume of a wing or area of a wing section \cite{MP2010}.
The displacement has to be maintained in order to guarantee that the optimized ship has the ability to transport the required payload.
Preserving the center of buoyancy serves two purposes: 
Firstly it supports maintaining the hydrostatic floating position and secondly it prevents the ship hull from being moved out of the computational domain.
Note that this is not sufficient to also account for changes of the floating position induced by the fluid dynamics.
This would also require to consider the rigid body dynamics, which balances the inertia and fluid dynamic forces and moments of the rigid vessel, and is ignored  herein.

The notation of this work employs $J'$ to indicate the derivative of a shape-based objective function $J$. The gradient of a function with respect to cartesian spatial coordinates is denoted by the nabla operator  $\nabla (\cdot) $ and the Jacobian refers to $D(\cdot)$.
The operator $A^*$ denotes the adjoint operator to $A$ and $\langle \cdot, \cdot \rangle$ is used to represent a dual pairing.

\par

\section{Mathematical Model}
We consider two and three dimensional flow domains $\Omega \subset \R^d$ ($d = 2,3$) which feature a boundary $\Gamma$ and obstacles $\Omega_{obs}$ with Lipschitz boundary $\Gamma_{obs}$ embedded in the flow domain.
We aim at minimizing the drag of the obstacle $\Omega_{obs}$ by deformation of a reference domain which at the same time serves as the initial configuration.
The domain is parameterized by a deformation field $u : \R^d \mapsto \R^d$.
We follow the approach in, e.g., \cite{SZ1992, DZ2011} where the domain is transformed by the perturbation of the identity.
Therewith the perturbed domain is defined by
\begin{equation}
	\tilde\Omega = (id + u)(\Omega) := \{x + u: x \in \Omega \}
\end{equation}
with $u \in W^{1,\infty}(\Omega, \R^d)$.
In the following we use the abbreviation $T(x) := x + u(x)$ and $T$ is a injective mapping with weakly differentiable inverse for sufficiently small $u$.
With this we approximate the shape derivative of  a shape function $J(\Omega)$ by the Fr{\'e}chet derivative of the mapping $W^{1,\infty}(\Omega,\R^d) \ni u \mapsto J((id + u)(\Omega))$, viz.  
\begin{equation}
	 J(\tilde{\Omega}) = J(\Omega) + J'(\Omega)u + o(\| u \|) \quad \text{for} \quad \| u \| \to 0
\end{equation}
where $u \mapsto J'(\Omega)u$ is linear regarding $u$.
To outline the central idea of this paper we first consider a generic constraint shape optimization problem
\begin{equation}
	\begin{aligned}
		\min_\Omega \quad \mathcal J(\Omega,y) \quad \text{s.t.} \quad e(\Omega,y) = 0\\
		\text{and} \quad g_i(\Omega,y) = 0, \quad i = 1, \dots, m.
	\end{aligned}
	\label{eq:GenericShapeOptProblem}
\end{equation}
where $e(\Omega, y)$ denotes the PDE constraint that describes the state $y$ and $g_i(\Omega, y), i = 1, \dots, m$ are a finite number of additional geometric constraints on the domain $\Omega$.
It is worth mentioning that the geometric constraints are restricting the shape of the domain $\Omega$ rather than contributing to the characterization of the state $y$. 
More precisely the state is fully described by the underlying boundary value problem $e(\Omega,y) = 0$ for a specific domain $\Omega$.
The change of the geometry, however, is restricted by the state and the solution of the PDE, respectively.
For the hydrodynamic problem at hand, this is described in greater detail in Section~\ref{sc:ComputationalModel}.
Assuming that the state constraint $e$ has a unique solution on $\Omega$ and thus the control-to-state map $\Omega \mapsto y(\Omega)$ exists, we obtain the reduced objective function $J(\Omega) := \mathcal J(\Omega, y(\Omega))$.
Furthermore we assume that the shape function $J(\Omega)$, as well as $e(\Omega,y)$ and $g_i(\Omega,y), \, i = 1, \dots, m$ are continuously Fr{\'e}chet differentiable.
Upon this we define the augmented Lagrange function
\begin{equation}
	L(\Omega, y, z, \lambda) = J(\Omega) + \langle z, e(\Omega, y) \rangle + \lambda^T g(\Omega, y)+ \frac{\tau}{2} \| g(\Omega, y) \|^2
	\label{eq:GenericLagrangeFunction}
\end{equation}
with the Lagrange multipliers $z$ and $\lambda = (\lambda_1, \dots, \lambda_m)^T \in \R^m$ and the penalty factor $\tau > 0$.
In the following $z$ is also referred to as the adjoint state.
\begin{algorithm}[htp]
	\caption{Shape optimization procedure}
	\label{alg:ShapeOptimizationProcedure}
	\begin{algorithmic}[1]
		\Repeat
		\State Compute $y$ at $u = 0$ with $e(\Omega, y) = 0$.
		\State Compute the adjoint state $z$ at $u = 0$  that fulfills
			\begin{equation}
				0 = \mathcal J_y(\Omega, y) \, \delta_y 
				+ \left\langle e_y(\Omega, y)^* \, z, \delta_y \right\rangle
				+ \lambda^T g_y (\Omega, y) \, \delta_y 
				+ \tau g(\Omega, y)^T g_y(\Omega,y) \delta_y 
				\qquad \forall \delta_y
				\label{eq:GenericAdjointEquation}
			\end{equation}
			\label{alg:line:GenericAdjointState}
		\State Find a descent direction $V$ at $u = 0$ by solving the minimization problem
			\begin{equation}
				\min_{V \in W^{1,\infty}, \, \| DV \| \leq 1} \; J'(\Omega)V = 
				\min_{V \in W^{1,\infty}, \, \| DV \| \leq 1} 
				\left\langle \mathcal J_u(\Omega, y) + e_u(\Omega, y)^* z - (\lambda - \tau g(\Omega, y) )^T g_u(\Omega, y), \, V \right\rangle \label{eq:GenericSteepestDescentDirection}
			\end{equation}
	           with an iterative scheme while successively updating $\lambda \gets \lambda + \tau \langle g_u(\Omega, y), V \rangle$ 
	           \label{alg:line:GenericDescentDirectionSubProblem}
	    \State Choose a sufficient step size $\epsilon > 0$
		\State Update the shape by applying the transformation
		    \begin{equation*}
		    	\Omega \gets (id + \epsilon V)(\Omega)
		    	\label{eq:GenericShapeUpdate}
			\end{equation*}
		\Until{converged}
	\end{algorithmic}
\end{algorithm}
Algorithm~\ref{alg:ShapeOptimizationProcedure} outlines the general shape optimization procedure which is based on an augmented Lagrange method of multipliers.
Usually the procedure would contain two nested loops, where the shape optimization problem by itself is solved several times with constant values of $\lambda$, and the update of the multiplier $\lambda$ is performed after each shape optimization loop. To reduce the related efforts in practical applications, the update of $\lambda$ is performed within the sub-optimization problem in \eqref{eq:GenericSteepestDescentDirection}, and  the previous value of $\lambda$ is used (as an approximation) in \eqref{eq:GenericAdjointEquation} for the current shape optimization step.
This allows to solve the shape optimization problem only once, but does not guarantee the convergence of the algorithm.
Indeed, all numerical experiments discussed in Section~\ref{sc:NumericalResults} show stable reductions of the objective functional while the geometric constraints are all met within a prescribed tolerance at each iteration.

The identification of a deformation field $V$  in accordance with Line~\ref{alg:line:GenericDescentDirectionSubProblem} of algorithm~\ref{alg:ShapeOptimizationProcedure} is  itself demanding.
To this end, we follow the approach suggested in \cite{DHH2022, MKS2021} and determine a descent direction by finding a minimizer of the $p$-Laplace relaxed problem
\begin{equation}
	\min_{V \in W^{1,p}} \frac{1}{p} \int_\Omega (\nabla V : \nabla V)^\frac{p}{2} \, d x
	+ \left\langle \mathcal J_u(\Omega, y) + e_u(\Omega, y)^* z - (\lambda - \tau g(\Omega, y) )^T g_u(\Omega, y), \, V \right\rangle \, . 
	\label{eq:GenericpLaplacianRelaxedDescentDirection}
\end{equation}
Relation (\ref{eq:GenericpLaplacianRelaxedDescentDirection}) approaches the limiting problem \eqref{eq:GenericSteepestDescentDirection} for $p \to \infty$, which characterizes the steepest descent direction in $W^{1,\infty}$-topology and therefore adheres to the first-order optimality condition
\begin{equation}
	\begin{aligned}
		&\int_\Omega (\nabla V : \nabla V)^\frac{p-2}{2} \, \nabla V : \nabla U \, d x \\
		&+ \left\langle \mathcal J_u(\Omega, y) + e_u(\Omega, y)^* z - (\lambda - \tau g(\Omega, y) )^T g_u(\Omega, y), \, U \right\rangle \geq 0 \quad \forall U \; . 
	\end{aligned}
	\label{eq:GenericpLaplacianRelaxedOptimalityConditions}
\end{equation}

\section{Computational Model}
\label{sc:ComputationalModel}

\begin{figure}[htp]
	\centering
	\includegraphics[page=1]{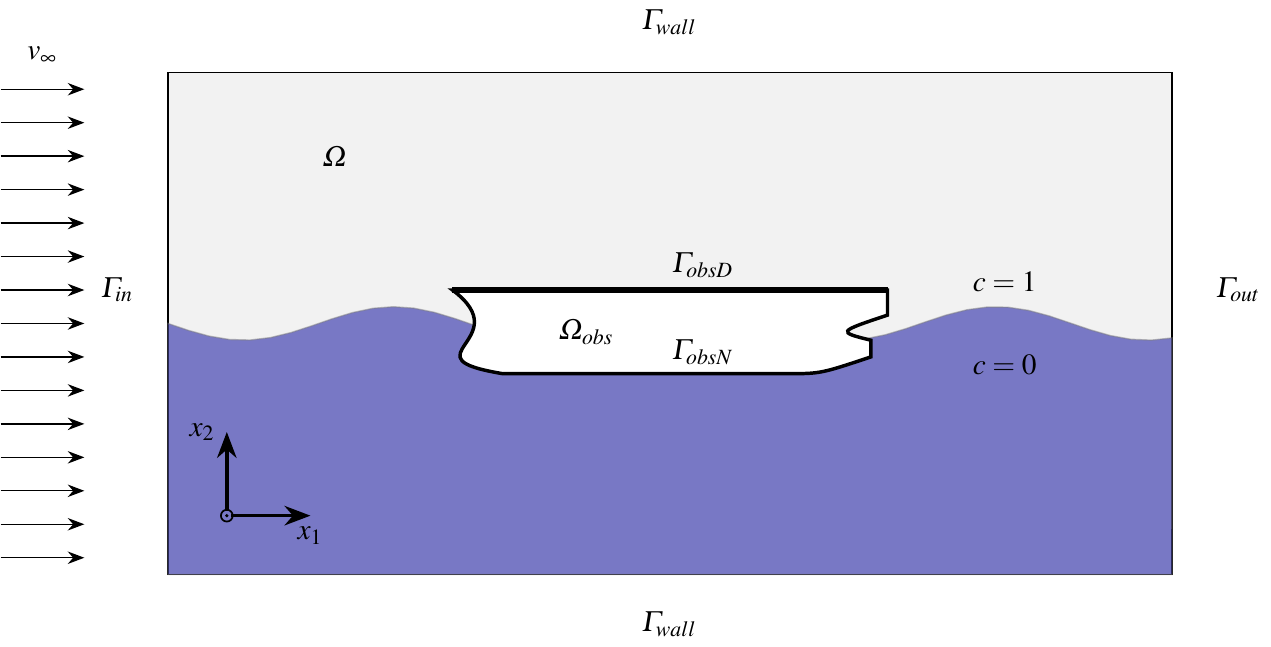}
	\caption{
		Illustration of the flow and obstacle domains and the boundaries. 
		The distribution of the fluid phase is indicated by the air volume concentration $c$, i.e. $c = 1[0]$ for air[water] filled regions.
	}
	\label{fig:Domain}
\end{figure}
As outlined in Figure~\ref{fig:Domain}, the flow domain is multiple-connected with one or several interior  boundaries in addition to a single outer boundary, and is occupied by two immiscible incompressible fluids, i.e. air and water.
The outer boundary is subdivided into the disjoint subsets inlet $\Gamma_{in}$, outlet $\Gamma_{out}$ and lateral as well as horizontal  walls $\Gamma_{wall}$.
The boundary $\Gamma_{obs}$ of the obstacle $\Omega_{obs}$ is fractioned into a nondeformable part $\Gamma_{obsD}$ and a deformable part $\Gamma_{obsN}$.
The optimization aims at minimal resistance of the obstacle $\Omega_{obs}$ by deforming the boundary $\Gamma_{obsN}$.
Because the flow domain $\Omega$ is the difference of the overall holdall domain and $\Omega_{obs}$ finding
a optimal shape of $\Omega_{obs}$ is equivalent to finding the optimal shape of the domain $\Omega$. 

On the one hand, the state is given by the velocity $v : \R_+ \times \Omega \to \R^d$, the total pressure $p : \R_+ \times \Omega \to \R$ and the indicator function/volume concentration $c : \R_+ \times \Omega \to [0,1]$ to distinguish between the two immiscible fluid phases air and water.
On the other hand the additional state variables kinetic turbulent energy $k : \R_+ \times \Omega \to \R_+$ and dissipation rate of kinetic turbulent energy $\omega : \R_+ \times \Omega \to \R$ are introduced for turbulence modeling.
Here we consider the Wilcox $k$-$\omega$ turbulence model \cite{W1998}.
The total pressure $p$ consists of the pressure $\mathfrak{p}$ and the hydrostatic pressure $- \rho \mathfrak{g} \cdot r$ with the acceleration due to gravity $\mathfrak{g}$ and the position $r$.
For the case sketched in Figure~\ref{fig:Domain} the gravitation is pointing in negative $x_2$-direction and thus with the basis vector of unit length $r = e_2$ the total pressure reads $p = \mathfrak{p} - \rho \mathfrak{g} x_2$.
For the turbulent two phase flow we consider the following RANS equations on $[0,T] \times \Omega$
\begin{equation}
	\begin{aligned}
		 \frac{\partial (\rho v)}{\partial t} + \div( \rho v \otimes v) - \div\left( (\mu + \mut) ( \nabla v + \nabla v^\T) \right) + \nabla p - f&= 0,\\
	    \div(v) &= 0,\\
	    \frac{\partial c}{\partial t}  + (v \cdot \nabla)c &= 0,\\
		\frac{\partial (\rho k)}{\partial t} + \div(\rho v k) - \div\left( (\mu + \sigma_k \mut) \nabla k \right) - P + \beta^\ast \rho k \omega &= 0,\\ 
		\frac{\partial (\rho \omega)}{\partial t} + \div(\rho v \omega) - \div\left( (\mu + \sigma_\omega \mut ) \nabla \omega \right) - \gamma \frac{\omega}{k} P + \beta \rho \omega^2 - \sigma_d \frac{\rho}{\omega} \nabla k \cdot \nabla \omega &= 0\\
	\end{aligned}
	\label{eq:PrimalURANSequations}
\end{equation}
where the turbulent production refers to $P = \mut (\nabla v + \nabla v^\T) : \nabla v$ and $f : \R_+ \times \Omega \to \R^d$ is a general source term that does not depend on the state.
The parameters $\mu, \rho > 0$ are the molecular viscosity and density of the respective fluids. 
They are composed from the bulk properties for air and water phase, i.e. $\rho_{\mathrm{air}}, \rho_{\mathrm{water}}$ and $\mu_{\mathrm{air}}, \mu_{\mathrm{water}}$, which are considered constant and the linear algebraic equation of state, viz. $\rho = \rho_{\mathrm{air}} c + \rho_{\mathrm{water}} (1-c)$ and $\mu = \mu_{\mathrm{air}} c + \mu_{\mathrm{water}} (1-c)$, holds.
This finally yields  a solenoidal velocity field with $\div(v)=0$.
The system is closed by the following set of initial and boundary conditions
\begin{equation}
	\begin{aligned}
		v = 0, \quad \frac{\partial c}{\partial n} = 0, \quad \frac{\partial k}{\partial n} = 0, \quad \frac{\partial \omega}{\partial n} &= 0 & \text{on} \, \Gamma_{obs} \; \text{and} \; \Gamma_{wall}, \\
		v = v_\infty, \quad c = c_\infty, \quad k = \frac{3}{2} \| v_\infty \|_2^2 \alpha^2, \quad \omega &= \frac{\rho k}{\mu \nu_t^+} & \text{on} \; \Gamma_{in}, \\
		(\mu + \mut) (\nabla v + \nabla v^\T) \cdot n = p n, \quad \frac{\partial c}{\partial n} = 0 \quad \frac{\partial k}{\partial n} = 0, \frac{\partial \omega}{\partial n} &= 0 & \text{on} \; \Gamma_{out}\\
		\text{and} \quad v(0) = v_0, \; c(0) = c_0, \; k(0) = \frac{3}{2} \| v_0 \|^2 \alpha^2, \; \omega(0) &= \frac{\rho k}{\mu \nu_t^+} & \text{in} \, \Omega.
	\end{aligned}
	\label{eq:PrimalBoundaryConditions}
\end{equation}
The parameters $\sigma_k, \sigma_\omega, \sigma_d, \gamma, \beta, \beta^\ast, \nu_t^+$ are real valued positive constants of the turbulence model and $\alpha \in (0,1]$.

The force vector $F$ acting on the boundary of the obstacle $\Gamma_{obs}$ is given by 
\begin{equation}
	F(\Gamma_{obs}) = \int_{\Gamma_{obs}} \left[ (\mu + \mut) \left( \nabla v + \nabla v^\T \right) \cdot n - p n \, \right] d s \, . 
	\label{eq:ObstacleForceSurface}
\end{equation}
The resistance or drag is associated with the component of the force in \eqref{eq:ObstacleForceSurface} in direction of the approaching flow which we assume to be aligned with the $x_1$-direction.
Hence, the resistance of the obstacle is given by the projection $- F(\Gamma_{obs}) \cdot e_1$ where $e_1$ is the basis vector of unit length in $x_1$-direction.
For the formulation of the shape derivative it is favorable to consider the volume formulation of the objective function.
Introducing a smooth extension $\eta : \Omega \to \R^d$ with $\eta |_{\Gamma_{obs}} \equiv - e_1$ and $\eta |_{\Gamma \setminus \Gamma_{obs}} \equiv 0$, one obtains the equivalent volume formulation through integration by parts of \eqref{eq:ObstacleForceSurface} \cite[Section 5.1]{BLUU2009}
\begin{equation}
	J(\Omega) = \int_\Omega \left[ (\rho (v \cdot \nabla)v - f) \cdot \eta + (\mu + \mut) \left( \nabla v + \nabla v^\T \right) : \nabla \eta - p \, \div(\eta) \right] d x \,.
	\label{eq:ObstacleForceVolume}
\end{equation}
We focus upon the steady state resistance and thus assume the flow to be stationary and all time derivatives in (\ref{eq:PrimalURANSequations}) vanish in the converged state.
In practice this means that the average over a sufficient pseudo-time/iteration period of the state variables is used to suppress minor remaining variations of the flow and the objective functional.
The geometric constraints for preserving the water displacement of $\Omega_{obs}$ and the center of buoyancy the obstacle are give by
\begin{equation}
	\begin{aligned}
		g_i(\Omega, c) &= \int_\Omega (1-c) x_i \; d x , \qquad i = 1, \dots, d \quad \text{and}\\
		g_{d+1}(\Omega, c) &= \int_\Omega (1- c) \; d x.
	\end{aligned}
	\label{eq:GeometricConstraints}
\end{equation}
The appearance of the concentration $c$ in (\ref{eq:GeometricConstraints}) secures that the displacement of the underwater hull of the vessel (water wetted part) is preserved rather than the volume of the whole hull. 
The formulation of the geometric constraints here differs substantially from previous work, e.g. \cite{OS2020} and \cite{MPRS2022}, as it depends on the geometry and the solution of a PDE.

Because the derivation as well as implementation of the adjoint problem corresponding to the primal problem \eqref{eq:PrimalURANSequations}~-~\eqref{eq:PrimalBoundaryConditions} holds several challenging aspects we follow common practice and neglecting the turbulence model for the adjoint system.
This simplification is also known as \textit{frozen turbulence assumption} \cite{DB2012, O2008, SR2013} 
where the state variables $k$ and $\omega$ are treated as constants when computing the derivative w.r.t. the state as well as the shape and the deformation field, respectively.
Hence, the state variable is considered $y = (v, p, c)$ in the following.

Hereon we define the augmented Lagrange function
\begin{equation}
	\begin{aligned}
		&L(\Omega, (v,p,c), (w, q, h), \gamma, \lambda) = \; J(\Omega) \\
		&\quad + \int_\Omega \left[ \left( \div\Big( \rho v \otimes v - (\mu + \rho \frac{k}{\omega}) ( \nabla v + \nabla v^\T) \Big) + \nabla p - f \right) \cdot w - \div(v) q + \div( v c ) h \right] \; d x \\
		&\quad + \int_{\Gamma_{obs}} \gamma \cdot v \; d x\\
		&\quad + \sum_{i=1}^d \lambda_i \int_\Omega (1-c) \, x_i  \; d x + \lambda_{d+1} \int_\Omega (1-c) \; d x\\
		&\quad + \frac{\tau}{2} \left( \sum_{i=1}^d \Big( \int_\Omega (1 - c) x_i \, d x \Big)^2 + \Big( \int_\Omega (1-c) \, d x \Big)^2 \right) \, , 
	\end{aligned}
	\label{eq:LagrangeFunction}
\end{equation}
where the multiplier $\gamma$ corresponds to the Dirichlet boundary conditions of the velocity $v = 0$ that hold on $\Gamma_{obs}$, and $\lambda = (\lambda_1, \dots, \lambda_d, \lambda_{d+1})$ is associated with the center of buoyancy and displacement constraint.
The adjoint state is characterized by the derivative of \eqref{eq:LagrangeFunction} w.r.t. the state $y = (v,p,c)$ which leads to the variational form
\begin{equation}
	\begin{aligned}
		0 =& \int_\Omega 
		\Big[ (\mu + \mut) (\nabla w + \nabla w^\T) : \nabla \delta_v
		- \rho v \cdot (\nabla w + \nabla w^\T) \cdot \delta_v\\
		&- \div(\delta_v) \, q
		- \delta_p \, \div(w) 
		+ \div(\delta_v c) h + \div(v \delta_c) \, h\\
		&+ \delta_c \Delta_\rho (v \otimes v) : \nabla w
		+ \delta_c \Delta_\mu (\nabla w + \nabla w^\T) : \nabla v \Big] \; d x\\
		&+ \int_{\Gamma_{obs}} \gamma \cdot \delta_v \; d x\\
		&+ \sum_{i=1}^d \lambda_i \int_\Omega - \delta_c \, x_i \; d x
		+ \lambda_{d+1} \int_\Omega - \delta_c \; d x\\
		&+ \tau \left( \sum_{i=1}^d \int_\Omega (1 - c) \, x_i \; d x \; \int_\Omega - \delta_c \, x_i \; d x \; + \; \lambda_{d+1} \int_\Omega (1 - c) \; d x \; \int_\Omega - \delta_c \; d x \right)\\
		& \forall \delta_y = (\delta_v, \delta_p, \delta_c)
	\end{aligned}
	\label{eq:AdjointEquation}
\end{equation}
where $\Delta_\rho = \rho_{\mathrm{air}} - \rho_{\mathrm{water}}$ and $\Delta_\mu = \mu_{\mathrm{air}} - \mu_{\mathrm{water}}$.
The boundary integrals vanish if the boundary conditions
\begin{equation}
	\begin{aligned}
		w &= -\eta & \text{on} \; \Gamma \setminus \Gamma_{out}\\
		\text{and} \quad (\mu + \mut) \big(\nabla w + \nabla w^\T \big) \cdot n &= (q - ch) n & \text{on} \; \Gamma_{out}.
	\end{aligned}
\end{equation}
hold and by choosing
\begin{equation}
	\gamma := - (\mu + \mut) \left(\nabla w + \nabla w^\T \right) \cdot n + (q - ch)n.
\end{equation}
In order to derive the directional derivative of the reduced objective $J'(\Omega)V$ we formally apply C{\'e}a's method.
For a detailed description see \cite[Section 4.6]{ADJ2020}.
In general the shape derivative of a objective function $J(\Omega)$ has a volume and an equivalent surface formulation.
For computational reasons it is favorable to consider the surface formulation even though it requires higher regularity of the solutions $(v,p,c)$ and $(w,q,h)$ of the primal and adjoint problem \eqref{eq:PrimalURANSequations}~-~\eqref{eq:PrimalBoundaryConditions} and \eqref{eq:AdjointEquation}, respectively.
Utilizing \cite[Theorem 4.2 and 4.3]{ADJ2020} and assuming that $(v, p, c)$ and $(w, q, h)$ have sufficient regularity we obtain
\begin{equation}
	\begin{aligned}
		J'(\Omega) V &= \int_{\Gamma_{obs}} \left( - (\mu + \mut) \frac{\partial w}{\partial n} \cdot \frac{\partial v}{\partial n} \right) \, V \cdot n \; d x\\
		&\quad + \sum_{i=1}^d (\lambda_i - \tau g_i(\Omega,c)) \int_{\Gamma_{obs}} (1 - c) x_i \, V \cdot n \; d x\\
		&\quad + (\lambda_{d+1} - \tau g_{d+1}(\Omega,c)) \int_{\Gamma_{obs}} (1 - c) \, V \cdot n \; d x.
	\end{aligned}
\end{equation}
As mentioned in Section~\ref{sc:Introduction} the deformation field is obtained from the directional shape derivative by solving the minimization problem \eqref{eq:GenericpLaplacianRelaxedDescentDirection}
\begin{equation}
	\begin{aligned}
		\min_{V \in W^{1,p}(\Omega, \R^d)}
		\frac{1}{p} \int_\Omega \left( DV : DV \right)^\frac{p}{2} \; d x + J'(\Omega)V.
	\end{aligned}
	\label{eq:DescentDirection}
\end{equation}
To ensure that the outer boundary remains unchanged the Dirichlet condition $u = 0$ holds almost everywhere on $\Gamma_{in} \cup \Gamma_{out} \cup \Gamma_{wall}$.
In addition parts of the obstacle may be fixed and thus $u = 0$ also holds a.e. on $\Gamma_{obsD}$ and natural boundary conditions hold on $\Gamma_{obsN}$ where the boundary is deformed.

To computing the shape deformation field $V$ characterized by the minimization problem \eqref{eq:DescentDirection} we suggest the procedure sketched in Algorithm~\ref{alg:PicardSolver}.
\begin{algorithm}[htp]
	\caption{Picard Iteration for Augmented $p$-Laplacian Problem}
	\label{alg:PicardSolver}
	\begin{algorithmic}[1]
		\State $\lambda \gets 0$, $p \gets 2$, $u \gets 0$
		\While{$p < p_{max}$}
		\State $k \gets 0$
		\Repeat
		\State Obtain a preliminary $\tilde V^k$ by solving the linearized problem
		\[\int_\Omega (\nabla V^{k-1} : \nabla V^{k-1})^\frac{p-2}{2} \, \nabla \tilde V^k : \nabla U \; d x + J'(\Omega)U \qquad \text{for all} \; U\]
		\label{alg:line:PicardLinearizedProblem}
		\State Relax update $V^k \gets V^k + \omega (\tilde V^k - V^{k-1})$ with $\omega \in (0,2)$
		\State Update multiplier $\lambda^k \gets \lambda^{k-1} + \tau \langle g_u(\Omega,y), V^k \rangle$
		\State $k \gets k + 1$
		\Until{$R^k = \| V^k - V^{k-1} \|_{L^2}^2 + \| \lambda_{bc}^k - \lambda_{bc}^{k-1} \|_2^2 + | \lambda_v^k - \lambda_v^{k-1}|^2 \leq tol$}\label{alg:line:PicardConverge}
		\State $p \gets p + p_{inc}$
		\EndWhile
	\end{algorithmic}
\end{algorithm}
For the sake of briefness we omit the precise solution method in Line~\ref{alg:line:PicardLinearizedProblem} of Algorithm~\ref{alg:PicardSolver} as the solution strategy does not depend on the discretization.
Nevertheless, in Section~\ref{sc:NumericalResults} we consider the finite volume method wherefore a formulation can be found by partial integration \cite{MKS2021}.

\section{Numerical Results}
\label{sc:NumericalResults}
Results presented in this paper are obtained from the finite volume procedure \textit{FreSCo\textsuperscript{+}} \cite{SR2011} for the KCS ship in model scale \cite{KCS2008} at Reynolds- and Froude numbers of $Re=1.43 \cdot 10^7$ and $Fn=0.26$.
Also the shape deformation field obtained from the $p$-Laplacian relaxed problem in \eqref{eq:GenericpLaplacianRelaxedOptimalityConditions} is approximated with a Picard iteration and finite volume discretization.

\begin{figure}[htp]
	\centering
	\includegraphics[width=0.7\columnwidth]{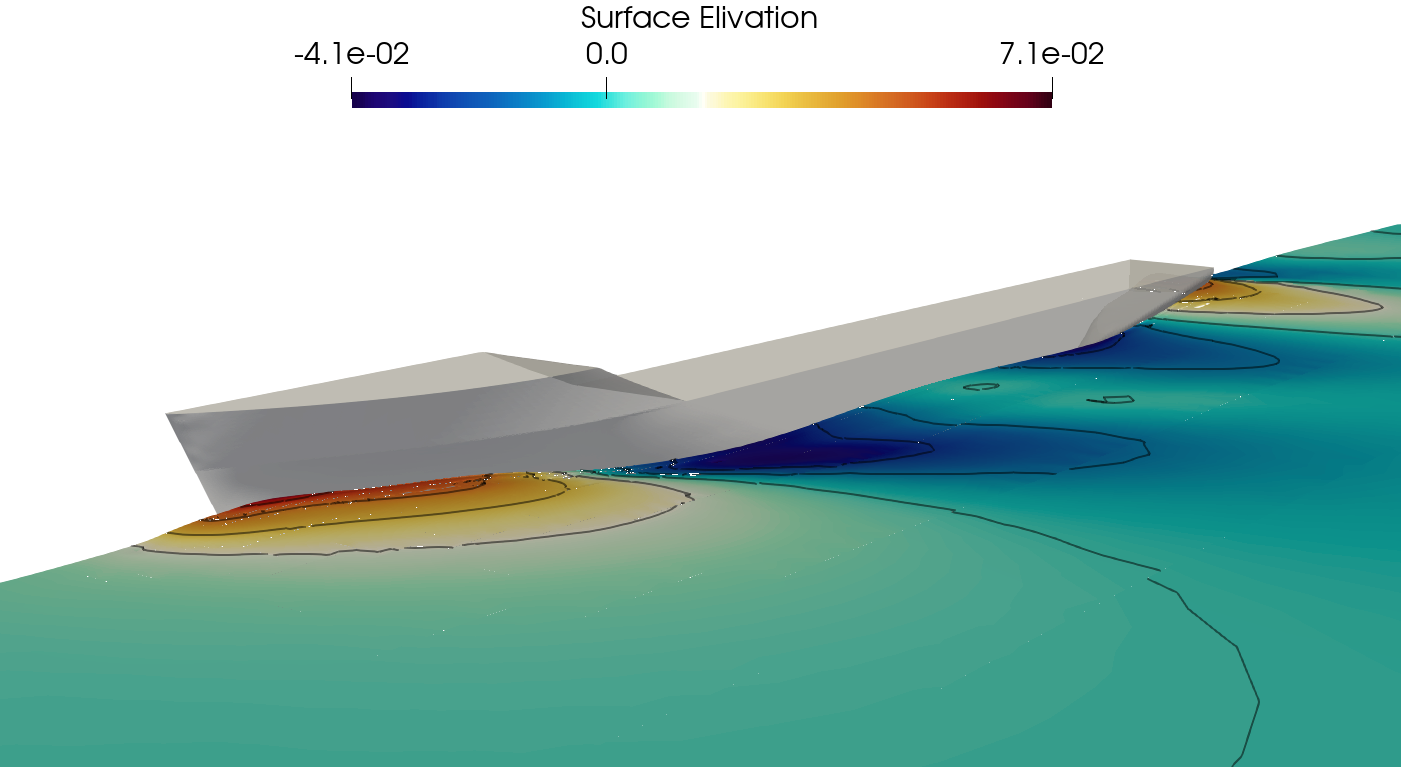}
	\caption{Initial hull shape and elevation of the free surface.}
	\label{fig:InitialFreeSurface}
\end{figure}

Figure~\ref{fig:InitialFreeSurface} shows the initial configuration with the hull of the KCS and free surface elevation.
We investigate two different cases to obtain the shape deformation field from \eqref{eq:DescentDirection} which differ in the boundary conditions along the hull.
Firstly, we consider the whole hull to be free for deformation, and secondly the air-wetted part of the hull remains fixed and only the underwater part $\Gamma_{obsN}$ of the hull is deformed.
In both studies the deck as well as the transom and a part of the propeller shaft remain fixed.
\begin{figure}[htp]
    \centering
    \includegraphics[page=1,width=0.7\linewidth]{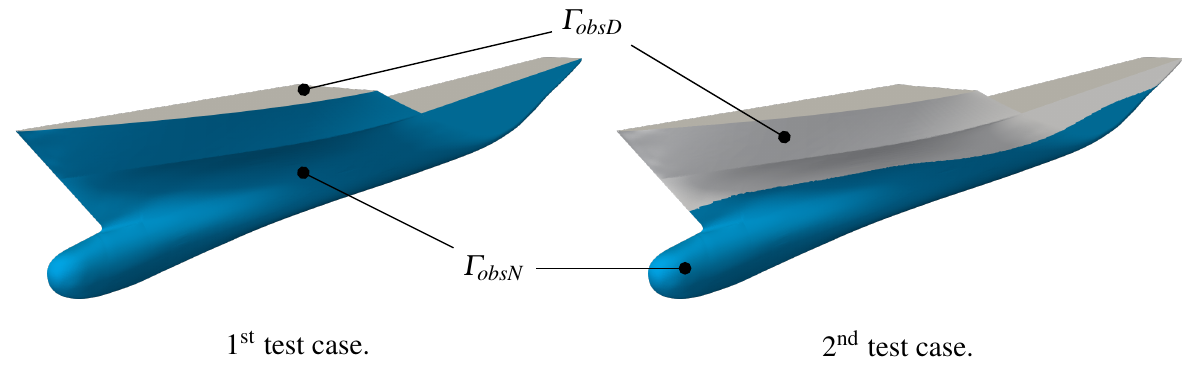}
    \caption{Boundary type layout for first test case (left) and second test case (right). The fixed part $\Gamma_{obsD}$ is colored in light gray and the deformed part $\Gamma_{obsN}$ in blue. The first test case assigns the complete hull to  $\Gamma_{obsN}$. The second test case only assigns the underwater hull to $\Gamma_{obsN}$, whereas air-wetted  hull parts above the water line belong to $\Gamma_{obsD}$.}
    \label{fig:BoundaryDeformationType}
\end{figure}
As stated in \cite{DHH2022, IL2005, MKS2021} the values for $p$ should be large in order to obtain a sufficient approximation for a descent direction in $W^{1,\infty}$.
However, due to the non-linearity of \eqref{eq:DescentDirection} the numerical computation for large values of $p$ is demanding  \cite{L2020} and we consider $p_{max} = 2.6$ as an upper value for both test cases. 
For Algorithm~\ref{alg:PicardSolver} to converge for $p > 2$ it requires a good initial guess $u^0$.
Therefore we consider the iteration over a sequence in $p=\{2, 2.3, 2.6\}$ to compute the initial guess for \eqref{eq:DescentDirection} with $p = p_{max}$ \cite{MKS2021, MPRS2022}.

\begin{figure}[htp]
	\centering
	\includegraphics[page=1,width=0.7\linewidth]{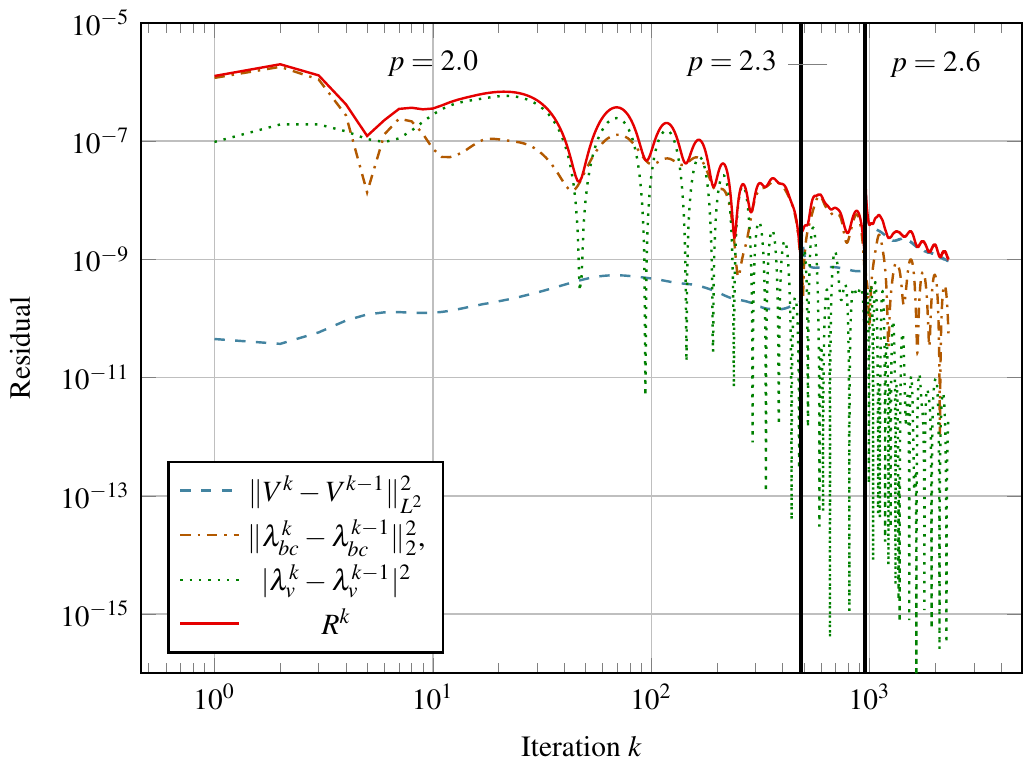}
	\caption{Residual plots of the $p$-Laplacian problem for the sequence of $p = \{2, 2.3, 2.6\}$. Displayed graphs refer to the residuals of the deformation field $V$, the Lagrangian multipliers $\lambda_{bc} = (\lambda_1, \dots, \lambda_d)$ and $\lambda_v = \lambda_{d+1}$ as well as the total residual $R$.}
	\label{fig:pLaplacianResidualPlot}
\end{figure}

To review Algorithm~\ref{alg:PicardSolver} we exemplary look at the first iteration of the shape optimization procedure in Algorithm~\ref{alg:ShapeOptimizationProcedure} for the first case.
Figure~\ref{fig:pLaplacianResidualPlot} shows the residuals of the procedure in Algorithm~\ref{alg:PicardSolver}.
The graphs display the individual contributions to the residual $R^k$ from Line~\ref{alg:line:PicardConverge} in Algorithm~\ref{alg:PicardSolver} for the tolerance $tol = 10^{-9}$.
The procedure is stable with the penalty factor $\tau = 10$ whereat the multipliers converge faster than $V^k \to V$.
The multiplier $\lambda$ and thus the choice of $\tau$ however heavily depends on the computed flow.

\begin{figure}[htp]
	\centering
	\includegraphics[page=1,width=0.7\linewidth]{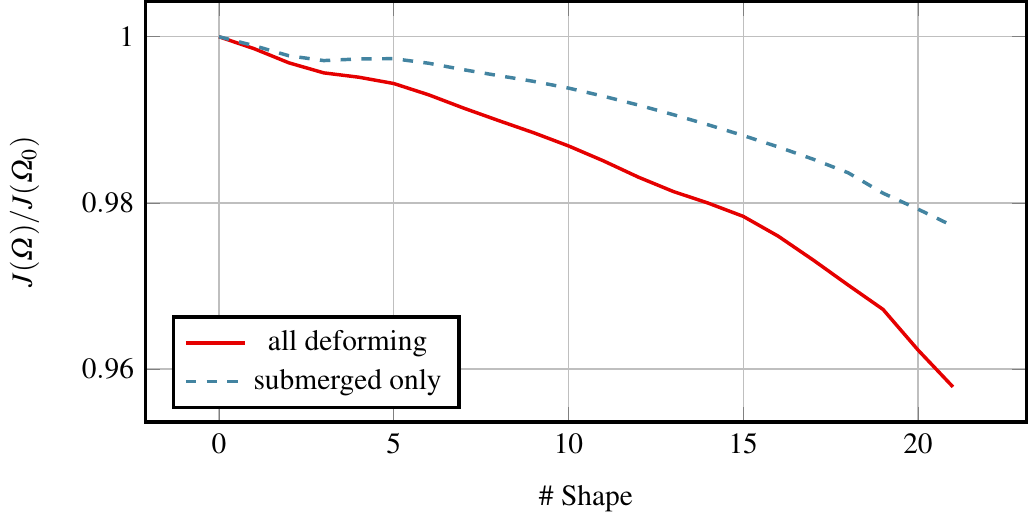}
	\caption{Evolution of the normalized drag force objective function obtained with $p = 2.6$ when deforming the whole hull (red line, test case 1) and only the underwater hull (blue line, test case 2).}
	\label{fig:FunctionalPlot}
\end{figure}

\begin{figure}[htp]
	\centering
	\includegraphics[page=1,width=0.7\linewidth]{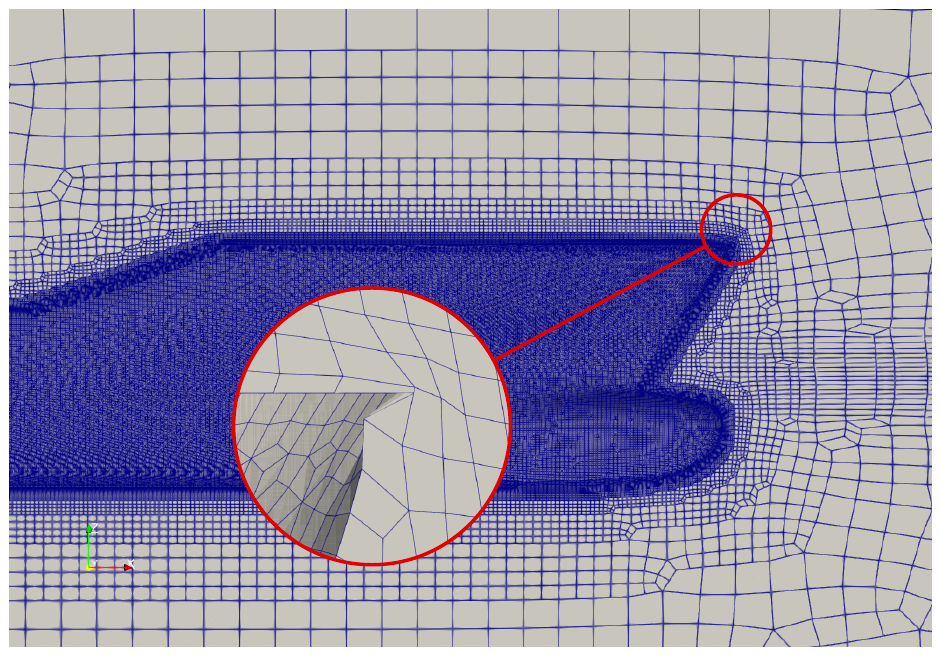}
	\caption{Illustration of the experienced grid deterioration in the bow region of the KCS container vessel after $22$ design iterations in combination with the first test case, where both the water and air wetted hull sections are deformed (cf. Fig. \ref{fig:BoundaryDeformationType}).}
	\label{fig:BowMesh}
\end{figure}

\begin{figure}
	\centering
	\includegraphics[page=1,width=0.8\columnwidth]{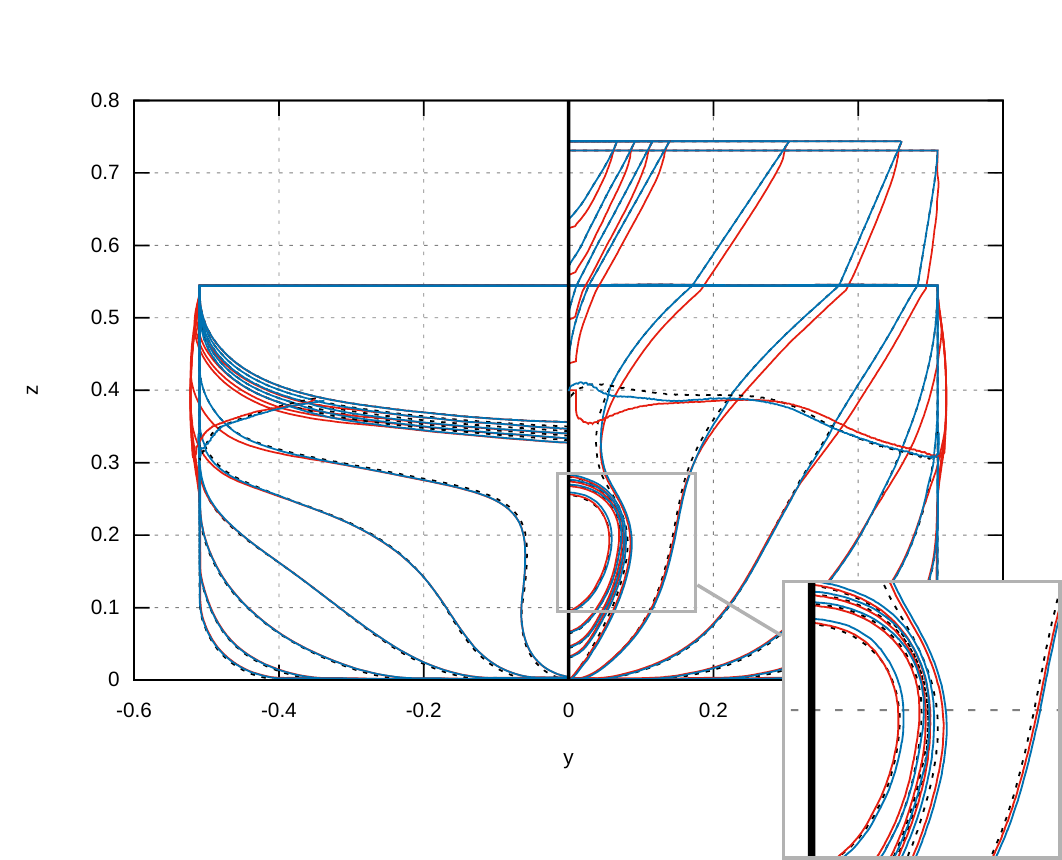}
	\caption{Body plan of the initial hull (dashed black line) and the optimized hulls of the 22nd design candidate predicted by the all deforming (red line) and only underwater hull deforming (blue line) approaches.}
	\label{fig:BodyPlanCompare}
\end{figure}

For both cases the normalized values of the objective function $- F(\Gamma_{obs}) \cdot e_1$ are shown in Figure~\ref{fig:FunctionalPlot}.
It can be observed that the functional values decline faster for the first test case where the deformation is not limited to the submerged part (solid line).
However, this is concomitant with large deformations at the intersection of the hull and the deck, particularly in the bow regime, cf. Fig. \ref{fig:BowMesh}. As outlined by the magnification in  Fig. \ref{fig:BowMesh}, we observe locally vanishing cell volumes after $22$ iterations in this regime and the simulations terminate. 
The issue only occurs for the first test case, where the whole vessel can deform. 
Using the second approach, the deformation is confined to the water-wetted surface and one could perform further iterations.
Figure~\ref{fig:BodyPlanCompare} compares the body plans of the initial (black) and the two modified designs of the 22\textsuperscript{nd} iteration. 
The two strategies predict virtually the same underwater hull deformations. However, differences occur when the free surface is approached, and more pronounced deformations are experienced in the first case, where the whole vessel can deform.   
Moreover, differences also occur in the bow regime, where the submerged only design (blue) predicts a stronger displacement in the upper part, as indicated by the magnification of the section lines close to the bow in Fig. ~\ref{fig:BodyPlanCompare}.

\section{Summary}
\label{sc:Summary}
We presented an algorithmic approach for fluid dynamic shape optimization of  floating ships exposed to turbulent two-phase flows under geometric constraints.
The main goal was to consider geometric constraints for displacement and the center of buoyancy in order to exclude trivial or undesirable optimal solutions in conjunction with shape updates which approximate the steepest descent direction in a Banach space.
The presented algorithm is based on the augmented Lagrange method of multipliers for the geometric constraints and the PDE constraints are treated utilizing the corresponding adjoint operator.
Numerical experiments were carried out based on the geometry of the KRISO Container Ship in model scale at realistic test conditions \cite{KCS2008}.

Results show that the suggested approach leads to deformation fields that fulfill the geometric constraints up to a predetermined tolerance.
However, the attainable drag reductions are limited by the degeneration of the computational grid and 
the discretization of the domain still becomes unfeasible coherent with the successive shape updates.
Because the domain of definition of the deformation field is not the holdall domain, but the flow domain (i.e. without the obstacle), it is globally not an injection and geometry overlaps can mathematically not be avoided.
Future research may thus consider a discretization of the entire domain including the interior of the obstacle.
Moreover, the algorithm may be applied to free-floating vessels subjected to rigid-body dynamics, which eventually gives a technically more realistic optimization problem.

\bibliographystyle{plain}
\bibliography{o2_references}

\begin{thebibliography}{10}

\bibitem{ADJ2020}
G.~Allaire, C.~Dapogny, and F.~Jouve.
\newblock Chapter 1 - shape and topology optimization.
\newblock In {\em Geometric Partial Differential Equations - Part II},
  volume~22 of {\em Handbook of Numerical Analysis}, pages 1 -- 132. Elsevier,
  2021.

\bibitem{BLUU2009}
C.~Brandenburg, F.~Lindemann, M.~Ulbrich, and S.~Ulbrich.
\newblock A continuous adjoint approach to shape optimization for navier stokes
  flow.
\newblock In {\em Optimal Control of Coupled Systems of Partial Differential
  Equations}, International Series of Numerical Mathematics. Birkhäuser,
  Basel, 2009.

\bibitem{DHH2022}
K.~Deckelnick, P.~J. Herbert, and M.~Hinze.
\newblock A novel $w^{1,\infty}y$ approach to shape optimization with lipschitz
  domains.
\newblock {\em ESAIM: COCV}, 28, 2022.

\bibitem{DZ2011}
M.~C. Delfour and J.-P. Zol\'esio.
\newblock {\em Shapes and Geometries}.
\newblock SIAM, 2011.

\bibitem{DB2012}
R.~P. Dwight and J.~Brezillon.
\newblock Effect of approximations of the discrete adjoint on gradient-based
  optimization.
\newblock {\em AIAA}, 44(12), 2012.

\bibitem{FP2020}
J.~H. Ferziger and M.~Peri{\'c}.
\newblock {\em Computational methods for fluid dynamics}.
\newblock Springer Berlin, Heidelberg, 4 edition, 2020.

\bibitem{HPUU2009}
M.~Hinze, R.~Pinnau, M.~Ulbrich, and S.~Ulbrich.
\newblock {\em Optimization with PDE Constraints}, volume~1.
\newblock 2009.

\bibitem{IL2005}
H.~Ishii and P.~Loreti.
\newblock Limits of solutions of p-laplace equations as p goes to infinity and
  related variational problems.
\newblock {\em SIAM journal on mathematical analysis}, 37(2):411 -- 437, 2005.

\bibitem{L2020}
S.~Loisel.
\newblock Efficient algorithms for solving the p-{L}aplacian in polynomial
  time.
\newblock {\em Numerische Mathematik}, 146(2):369--400, 2020.

\bibitem{KCS2008}
Maritime and Ocean Engineering~Research Institute.
\newblock Kriso container ship geometry.
\newblock \url{http://www.simman2008.dk/KCS/kcs_geometry.htm}.

\bibitem{MP2010}
B.~Mohammadi and O.~Pironneau.
\newblock {\em Applied shape optimization for fluids}, volume~2.
\newblock Oxford University Press, 2010.

\bibitem{MKS2021}
P.~M. M\"uller, N.~K\"uhl, M.~Siebenborn, K.~Deckelnick, M.~Hinze, and T.~Rung.
\newblock A novel p-harmonic descent approach applied to fluid dynamic shape
  optimization.
\newblock {\em Struct Multidic Optim}, 64, 2021.

\bibitem{MPRS2022}
P.~M. M{\"u}ller, J.~Pinz{\'o}n, Thomas Rung, and Martin Siebenborn.
\newblock A scalable algorithm for shape optimization with geometric
  constraints in banach spaces.
\newblock \url{https://arxiv.org/abs/2205.01912}, 2022.

\bibitem{OS2020}
S.~Onyshkevych and M.~Siebenborn.
\newblock Mesh quality preserving shape optimization using nonlinear extension
  operators.
\newblock {\em Journal of Optimization Theory and Applications}, 189:291--316,
  2020.

\bibitem{O2008}
C.~Othmer.
\newblock A continuous adjoint formulation for the computation of topological
  and surface sensitivities of ducted flows.
\newblock {\em Numerical Methods in Fluids}, 58(8):861--877, 2008.

\bibitem{OG2009}
E.~{\"O}zkaya and N.~R. Gauger.
\newblock Single-step one-shot aerodynamic shape optimization.
\newblock In {\em Optimal Control of Coupled Systems of Partial Differential
  Equations}, pages 191--204. Birkh{\"a}user Basel, 2009.

\bibitem{SZ1992}
J.~Sokolovski and J.-P. Zol\'esio.
\newblock {\em Introduction to shape optimization}.
\newblock Springer-Verlag, 1992.

\bibitem{SR2011}
A.~St\"uck and T.~Rung.
\newblock Adjoint rans with filtered shape derivatives for hydrodynamic
  optimisation.
\newblock {\em Computers \& Fluids}, 47(1):22--32, 2011.

\bibitem{SR2013}
A.~St{\"u}ck and T.~Rung.
\newblock Adjoint complement to viscous finite-volume pressure-correction
  methods.
\newblock {\em Journal of Computational Physics}, 248:402--419, 2013.

\bibitem{W1998}
David~C Wilcox et~al.
\newblock {\em Turbulence modeling for CFD}, volume~2.
\newblock DCW industries La Canada, CA, 1998.

\end{thebibliography}

\paragraph{Acknowledgements}
\label{sc:Acknowledgements}
The authors acknowledge the support by the Deutsche Forschungsgemeinschaft (DFG) within the Research Training Group GRK 2583 “Modeling, Simulation and Optimization of Fluid Dynamic Applications”.
The computations were performed with resources provided by the North-German Super-computing Alliance (HLRN).

\paragraph{Replication of results}
The geometry of the KCS is available at \url{http://www.simman2008.dk/KCS/kcs_geometry.htm}. 
A proprietary software is used for mesh generation.
Computations are carried out with the in-house finite volume code \textit{FreSCo\textsuperscript{+}}.

\end{document}